\documentclass[11pt,fleqn]{article}
\usepackage{amsfonts, epsfig, amssymb, amsmath}
\usepackage[normalem]{ulem}
\usepackage{color}
\newcommand{\ol}{\setlength{\itemsep}{0pt.}\begin{enumerate}}
\newcommand{\eol}{\end{enumerate}\setlength{\itemsep}{-\parsep}}
\newcommand{\ignore}[1]{}
\title{Eigenvalues and eigenfunctions of a Hamming ball}
\author{Amit Avni\thanks{School of Engineering and Computer Science,
The Hebrew University of Jerusalem,
Jerusalem 91904, Israel.
Research partially supported by ISF
grant 921/22. } ~and Alex Samorodnitsky\thanks{School of Engineering and Computer Science,
The Hebrew University of Jerusalem,
Jerusalem 91904, Israel.
Research partially supported by ISF
grant 921/22. }
}

\begin{document}
\date{}
\maketitle


\newtheorem{THEOREM}{Theorem}[section]
\newenvironment{theorem}{\begin{THEOREM} \hspace{-.85em} {\bf :}
}%
                        {\end{THEOREM}}
\newtheorem{LEMMA}[THEOREM]{Lemma}
\newenvironment{lemma}{\begin{LEMMA} \hspace{-.85em} {\bf :} }%
                      {\end{LEMMA}}
\newtheorem{COROLLARY}[THEOREM]{Corollary}
\newenvironment{corollary}{\begin{COROLLARY} \hspace{-.85em} {\bf
:} }%
                          {\end{COROLLARY}}
\newtheorem{PROPOSITION}[THEOREM]{Proposition}
\newenvironment{proposition}{\begin{PROPOSITION} \hspace{-.85em}
{\bf :} }%
                            {\end{PROPOSITION}}
\newtheorem{DEFINITION}[THEOREM]{Definition}
\newenvironment{definition}{\begin{DEFINITION} \hspace{-.85em} {\bf
:} \rm}%
                            {\end{DEFINITION}}
\newtheorem{EXAMPLE}[THEOREM]{Example}
\newenvironment{example}{\begin{EXAMPLE} \hspace{-.85em} {\bf :}
\rm}%
                            {\end{EXAMPLE}}
\newtheorem{CONJECTURE}[THEOREM]{Conjecture}
\newenvironment{conjecture}{\begin{CONJECTURE} \hspace{-.85em}
{\bf :} \rm}%
                            {\end{CONJECTURE}}
\newtheorem{MAINCONJECTURE}[THEOREM]{Main Conjecture}
\newenvironment{mainconjecture}{\begin{MAINCONJECTURE} \hspace{-.85em}
{\bf :} \rm}%
                            {\end{MAINCONJECTURE}}
\newtheorem{PROBLEM}[THEOREM]{Problem}
\newenvironment{problem}{\begin{PROBLEM} \hspace{-.85em} {\bf :}
\rm}%
                            {\end{PROBLEM}}
\newtheorem{QUESTION}[THEOREM]{Question}
\newenvironment{question}{\begin{QUESTION} \hspace{-.85em} {\bf :}
\rm}%
                            {\end{QUESTION}}
\newtheorem{REMARK}[THEOREM]{Remark}
\newenvironment{remark}{\begin{REMARK} \hspace{-.85em} {\bf :}
\rm}%
                            {\end{REMARK}}

\newcommand{\thm}{\begin{theorem}}
\newcommand{\lem}{\begin{lemma}}
\newcommand{\pro}{\begin{proposition}}
\newcommand{\dfn}{\begin{definition}}
\newcommand{\rem}{\begin{remark}}
\newcommand{\xam}{\begin{example}}
\newcommand{\cnj}{\begin{conjecture}}
\newcommand{\mcnj}{\begin{mainconjecture}}
\newcommand{\prb}{\begin{problem}}
\newcommand{\que}{\begin{question}}
\newcommand{\cor}{\begin{corollary}}
\newcommand{\prf}{\noindent{\bf Proof:} }
\newcommand{\ethm}{\end{theorem}}
\newcommand{\elem}{\end{lemma}}
\newcommand{\epro}{\end{proposition}}
\newcommand{\edfn}{\bbox\end{definition}}
\newcommand{\erem}{\bbox\end{remark}}
\newcommand{\exam}{\bbox\end{example}}
\newcommand{\ecnj}{\bbox\end{conjecture}}
\newcommand{\emcnj}{\bbox\end{mainconjecture}}
\newcommand{\eprb}{\bbox\end{problem}}
\newcommand{\eque}{\bbox\end{question}}
\newcommand{\ecor}{\end{corollary}}
\newcommand{\eprf}{\bbox}
\newcommand{\beqn}{\begin{equation}}
\newcommand{\eeqn}{\end{equation}}
\newcommand{\wbox}{\mbox{$\sqcap$\llap{$\sqcup$}}}
\newcommand{\bbox}{\vrule height7pt width4pt depth1pt}
\newcommand{\qed}{\bbox}
\def\sup{^}

\def\H{\{0,1\}^n}

\def\S{S(n,w)}

\def\g{g_{\ast}}
\def\t{t^{\ast}}
\def\xop{x_{\ast}}
\def\y{y_{\ast}}
\def\z{z_{\ast}}

\def\f{\tilde f}

\def\n{\lfloor \frac n2 \rfloor}

\def \E{\mathop{{}\mathbb E}}
\def \R{\mathbb R}
\def \Z{\mathbb Z}
\def \F{\mathbb F}
\def \T{\mathbb T}

\def \x{\textcolor{red}{x}}
\def \r{\textcolor{red}{r}}
\def \Rc{\textcolor{red}{R}}

\def \noi{{\noindent}}

\def \iff{~~~~\Longleftrightarrow~~~~}

\def\myblt{\noi --\, }

\def \queq {\quad = \quad}

\def\<{\left<}
\def\>{\right>}
\def \({\left(}
\def \){\right)}

\def \e{\epsilon}
\def \l{\lambda}

\def\Tp{Tchebyshef polynomial}
\def\Tps{TchebysDeto be the maximafine $A(n,d)$ l size of a code with distance $d$hef polynomials}
\newcommand{\rarrow}{\rightarrow}

\newcommand{\larrow}{\leftarrow}

\overfullrule=0pt
\def\setof#1{\lbrace #1 \rbrace}

\begin{abstract}
We describe the eigenvalues and the eigenspaces of the adjacency matrices of subgraphs of the Hamming cube induced by Hamming balls, and more generally, by a union of adjacent concentric Hamming spheres. 

As a corollary, we extend the range of cardinalities of subsets of the Hamming cube for which Hamming balls have essentially the largest maximal eigenvalue (among all subsets of the same size). We show that this holds even when the sets in question are large, with cardinality which is an arbitrary subconstant fraction of the whole cube. 
\end{abstract}

\section{Introduction}

\noi This paper studies the eigenvalues and the eigenfunctions of subgraphs of the Hamming cube induced by Hamming balls or, more generally, by a union of adjacent concentric Hamming spheres. The maximal eigenvalue of the Hamming ball of a given radius is an interesting object of study. It plays an important role in the proof of the first linear programming bound for linear binary codes, due to \cite{ft}. Furthermore, \cite{ft, BLL, sam-log-sob} consider the question of whether the Hamming ball has the largest (or essentially largest) maximal eigenvalue among all the subsets of the Hamming cube with the same cardinality (at least for most values of this cardinality). This has been partially validated in \cite{BLL, sam-log-sob}.

\noi The role played by the maximal eigenvalue of the Hamming ball in the argument of \cite{ft} is similar to the role played by the first root of the appropriate Krawtchouk polynomial in the original proof of the first linear programming bound for binary codes (\cite{mrrw}). The estimates in \cite{ft} show these two values to be closely related. One of the starting goals in this paper was to make the connection between these two notions more precise. It has been observed in \cite{llc} that if $\l$ is the maximal eigenvalue of a Hamming ball of radius $r$, then $\frac{n-\l}{2}$ is a root of a Krawtchouk polynomial. Here we will show (see Corollary~\ref{cor:max eigenvalue} below for a precise statement) that $\frac{n-\l}{2}$ is in fact the first root of the Krawtchouk polynomial $K^{(n)}_{r+1}$  (all notions will be defined below).

\noi Furthermore, we describe all the eigenvalues of the subgraphs induced by Hamming balls or by a union of concentric Hamming spheres in terms of roots of Krawtchouk polynomials. We also provide a description of the corresponding eigenfunctions. At the intuitive level, since these subgraphs are induced by a collection of spheres which are all embedded in the same Hamming cube, they should inherit some of the properties of both the Hamming cube and the Hamming sphere. In a sense, we provide a certain 'algebraic' validation of this intuition. Recall that the Hamming cube and the Hamming sphere are distance-transitive spaces, and hence support an association scheme structure, the Hamming association scheme for the cube and the Johnson association scheme for the sphere. We describe the eigenvalues of the subgraphs under consideration in terms of {\it Krawtchouk polynomials}, namely the {\it $Q$-polynomials} of the Hamming scheme, and (speaking somewhat imprecisely) the eigenfunctions in terms of the {\it Hahn polynomials}, namely the $Q$-polynomials of the Johnson scheme.

\noi We start with a brief description of the notions involved, following \cite{LW, MT, mrrw}.

\subsubsection{Doubly transitive spaces and association schemes}
\label{subsubsec:doubly transitive}

\noi We consider two metric spaces, the Hamming cube $\H$ endowed with the Hamming distance, $|x-y| = |\{j: x_j \not = y_j\}|$, and the $n$-dimensional Hamming sphere $S(n,i)$ of radius $i \le n/2$, which is a subset of $\H$ consisting of all points at distance $i$ from $0$. Both of these spaces are {\it distance transitive}, namely they are acted on by a group of isometries $I$, so that for any two pairs of points $\(x_1,y_1\)$ and $\(x_2,y_2\)$ in the space with $|x_1 - y_1| = |x_2 - y_2|$, there is an isometry $\alpha \in I$ taking $x_1$ to $x_2$ and $y_1$ to $y_2$. In general, let $(X,d)$ be a finite distance-transitive metric space, and let $I$ be its isometry group. Each isometry $\alpha \in I$ corresponds to a linear operator $T_{\alpha}$ acting on the functions on $X$, which is defined in the following way: $\(T_{\alpha} f\)(x) = f(\alpha(x))$. The linear space of operators commuting with all the isometries $\{T_{\alpha}\}_{\alpha \in I}$ is closed under composition. Furthermore, all these operators commute. This space is called the {\it Bose-Mesner algebra} of $X$, or of the related association scheme. The dimension of this space is the number $k$ of distinct distances in $X$. Let these distances be $0 = d_0, d_1,...,d_{k-1}$. The Bose-Mesner algebra has two important bases. The operators in the first basis are represented in the standard basis of the functions on $X$ by the {\it adjacency matrices} $A_0,...,A_{k-1}$, where $A_j$ is an $|X| \times |X|$ matrix, given by $A_j\(x_1,x_2\) = \left\{\begin{array}{cc} 1 & d\(x_1,x_2\) = d_j \\ 0 & \mathrm{otherwise} \end{array}\right.$. The second basis is formed by the orthogonal projections $E_0,...,E_t$ on eigenspaces $V_0,V_1,...,V_t$ of the space $X$. These are the joint eigenspaces of the adjacency matrices of $X$, and we have $V_0 \oplus V_1 \cdots \oplus V_t = \R^X$. The transformation matrices between these two bases of the Bose-Mesner algebra are called the $P$ and $Q$-eigenmatrices of the scheme.

\noi We also need the notions of a {\it spherical} and a {\it zonal spherical} function.

\dfn
\label{dfn:spherical}

Let $(X,d)$ be a finite distance-transitive metric space and let $o \in X$. A function $f$ on $X$ is called spherical around $o$ if its value at a point $x \in X$ depends only on $d(x,o)$. A function $f \not = 0$ is the $j^{th}$ zonal spherical function if it is spherical and it lies in the eigenspace $V_j$.
\edfn

\noi We have the following well-known fact.

\lem
\label{lem:spherical zonal}
For each $o \in X$, and for each $0 \le j \le k-1$, there is a unique, up to a multiplicative scalar, zonal spherical function $F_{o,j}$ around $o$ in $V_j$. Moreover, the functions $\{F_{o,j}\}_{o \in X}$ span $V_j$.
\elem

\noi Specializing to the Hamming cube and the Hamming sphere, we have the following.

\begin{itemize}

\item {\it The Hamming cube}. The isometries of $\H$ are generated by the permutations in $S_n$, with $\sigma: \(x_1 \ldots x_n\) \rarrow \(x_{\sigma(1)} \ldots x_{\sigma(n)}\)$ and by addition in the group $\F_2^n$, with $a: x \rarrow x + a$. The distinct distances in $\H$ are $0,1,...,n$ and hence the Bose-Mesner algebra is $(n+1)$-dimensional. To describe the eigenspaces $V_0,...,V_n$, we recall some basic notions in the Fourier analysis on the boolean cube (see \cite{O'Donnel}). For $\gamma \in \H$, define the Walsh-Fourier character $W_{\gamma}$ on $\H$ by setting $W_{\gamma}(y) = (-1)^{\sum \gamma_i y_i}$, for all $y \in \H$. The {\it weight} of the character $W_{\gamma}$ is the Hamming weight $|\gamma| = |\{j: \gamma_j \not = 0\}|$ of $\gamma$.  The characters $\{W_{\gamma}\}_{\gamma \in \H}$ form an orthonormal basis in the space of real-valued functions on $\H$, under the inner product $\<f, g\> = \frac{1}{2^n} \sum_{x \in \H} f(x) g(x)$. We can now describe the eigenspaces of $\H$. We have $V_j = \<W_{\gamma}\>_{|\gamma| = j}$, $j=0,...,n$. Note that $\dim\(V_j\) = {n \choose j}$.

    The zonal spherical functions on $\H$ can be described as follows. For $0 \le j \le n$, let $F_j$ be the sum of all Walsh-Fourier characters of weight $j$, that is $F_j = \sum_{|\gamma| = j} W_{\gamma}$. It is easy to see that $F_j(x)$ depends only on the Hamming weight $|x|$ of $x$, and hence it is a $j^{th}$ zonal spherical function around $0$. Viewed as a univariate function on the integer points $0,...,n$, the function $F_j$ is given by the restriction to $\{0,...,n\}$ of the univariate polynomial $K_j = \sum_{\ell=0}^j (-1)^{\ell} {x \choose {\ell}} {{n-x} \choose {j-\ell}}$ of degree $j$. That is, $F_j(x) = K_j(|x|)$. The polynomial $K_j$ is the $j^{th}$ {\it Krawtchouk polynomial}. Abusing notation, we will also call $F_j$ the $j^{th}$ Krawtchouk polynomial, and write $K_j$ for $F_j$ when the context is clear. Since we will consider Hamming cubes $\{0,1\}^m$ of varying dimension, we will add a superscript denoting this dimension, if required. Thus $\{K^{(m)}_j\}_{j=0}^m$ are the $m+1$ zonal spherical functions around $0$ in $\{0,1\}^m$. An equivalent statement in terms of association schemes is that the Krawtchouk polynomials are the $Q$-polynomials of the Hamming association scheme (they describe the entries in the $Q$-eigenmatrix of that scheme).

\item {\it The Hamming sphere $S(n,i)$}. The isometries of $S(n,i)$ are given by the permutations in $S_n$, with $\sigma: \(x_1,...,x_n\) \rarrow \(x_{\sigma(1)}...x_{\sigma(n)}\)$. The distinct distances in $S(n,i)$ are $0,2,...,2i$ and hence the Bose-Mesner algebra is $(i+1)$-dimensional. The eigenspaces of $S(n,i)$ can be described as follows (see Example 30.7 in \cite{LW}). First, we define some functions which will play an important role below.

    \dfn
    \label{dfn:functions g}
    Let $z \in \H$, $|z| \le i$. We define $g_z$ as a function on $S(n,i)$ given by: $g_z(x) = \left\{\begin{array}{ccc} 1 & \mathrm{if} & z \subseteq x \\ 0 & \mathrm{otherwise} \end{array} \right.$
    \edfn

    For $0 \le j \le i$, let $U_j$ be the span of the functions $g_z$, $|z| \le j$. Clearly $U_0 \subseteq U_1 \subseteq ... \subseteq U_i = \R^{S(n,i)}$. Furthermore, it can be shown that $\dim\(U_j\) = {n \choose j}$. Define $U_0 = V_0$ to be the space of constant functions on $S(n,i)$, and let $V_j = U_j \cap \(U_{j-1}\)^{\perp}$ for $1 \le j \le i$. Then $V_0,...,V_i$ are the eigenspaces of the scheme. Note that $\dim\(V_j\) = {n \choose j} - {n \choose {j-1}}$.

    The zonal spherical functions on $S(n,i)$ can be described as follows. Let $o \in S(n,i)$. Then the $j^{th}$ zonal spherical function around $o$ is given by the restriction to the integer points $0,...,i$ of the {\it Hahn polynomial} $Q_j(x)$ of degree $j$. Hahn polynomials $Q_0,...,Q_i$ are the {\it $Q$-polynomials} of the Johnson association scheme. (See e.g., \cite{mrrw}).

    \end{itemize}

\subsection{Our results}

\noi We start with introducing the notion of {\it semi-symmetric functions} on the Hamming cube and the Hamming sphere.

\dfn
\label{dfn:semi-symmetric}

\begin{itemize}

\item A function $f$ on $\H$ is semi-symmetric around $y \in \H$ if its value at a point $x$ depends only on the Hamming distances $|x| = |x-0|$ and $|x-y|$. An important special case is $y = 0$. In this case, $f(x)$ depends only on $|x|$, and it is a spherical function around $0$.

\item A function $f$ on $S(n,i)$ is semi-symmetric around $y \in \H$ if its value at a point $x$ depends only on the Hamming distance $|x-y|$. An important special case is $y \in S(n,i)$. In this case, $f$ is a spherical function around $y$.

\end{itemize}

\edfn

\noi  {\bf Notation}: We will denote by $S_y$ the space of semi-symmetric functions around $y$ in $\H$ and by $S_{y,i}$ the space of semi-symmetric functions around $y$ in $S(n,i)$. Note that $S_{y,i}$ is obtained by restricting the functions in $S_y$ to $S(n,i)$.

\noi The following simple claim will be used several times below.

\lem
\label{lem:y-invariant}
Both $S_y$ and $S_{y,i}$  are linear subspaces of functions on $\H$ and $S(n,i)$ respectively, containing specifically the functions which are invariant under permutations in $S_n$ which preserve the point $y$.
\elem

\prf
\eprf

\noi We prove the following claim, in analogy with Lemma~\ref{lem:spherical zonal}.

\pro
\label{pro:semi-symmetric}

Let $0 \le j \le i \le n/2$. Let $V_j$ be an eigenspace of $S(n,i)$. Then, up to a multiplicative scalar, for each $y \in \H$, $|y| = j$, there is a unique non-zero semi-symmetric function $F_{y,j}$ around $y$ in $V_j$. We will refer to such functions as "zonal semi-symmetric functions". Furthermore, the functions $\{F_{y,j}\}_{y \in S(n,j)}$ span $V_j$. 
\epro

\ignore{
\noi Note that the notion of a zonal semi-symmetric function generalizes that of a spherical zonal function, and Proposition~\ref{pro:semi-symmetric} generalizes Lemma~\ref{lem:spherical zonal}.
}

\noi We proceed to study the spectral properties of subgraphs of $\H$ induced by Hamming balls, starting with some preliminaries. 

\noi Recall that $\H$ can be viewed as a graph in which two vertices are connected by an edge iff their Hamming distance is $1$. Let $0 \le r \le n/2$. The $n$-dimensional Hamming ball of radius $r$ around $0$ is a subset of $\H$ defined as $B = B(n,r) = \cup_{i=0}^r S(n,i)$. That is, $B$ is a union of $r+1$ concentric Hamming spheres. Let $A = A(n,r)$ be the adjacency matrix of the subgraph of $\H$ induced by the vertices in $B$. We will describe the eigenvalues of $A$ in terms of roots of Krawtchouk polynomials. Recall that we denote by $K^{(m)}_k$ the $k^{th}$ Krawtchouk polynomial on $\{0,1\}^m$. Krawtchouk polynomials $\{K^{(m)}_k\}_{k=0}^m$ form a family of orthogonal polynomials (with respect to the binomial measure $w_i = \frac{{m \choose i}}{2^m}$ on $0,1,...,m$) and hence (\cite{Szego}) all of their roots are real and distinct and lie in the interval $(0,m)$. Let ${\cal R}(m,k)$ be the set of roots of $K^{(m)}_k$. 

\noi {\bf Notation}. For $0 \le t \le r$, let 
\[
\Lambda_t ~=~ 2 \cdot {\cal R}(n-2t,r-t+1) - (n-2t),
\]
where for $X \subseteq \R$ and $s, m \in \R$, we write $m X + s$ for the set $\{mx + s\}_{x \in X}$.

\noi Then we have the following main claim. 

\thm
\label{thm:main}

\begin{enumerate}

\item  For each choice of $t = 0,...,r$ and for each of the $r-t+1$ distinct real numbers $\l \in \Lambda_t$, the matrix $A$ has an eigenspace $W_{t}(\l)$ of dimension ${n \choose t} - {n \choose {t-1}}$ corresponding to the eigenvalue $\l$.

\item A real number $\l$ is an eigenvalue of $A$ if and only if it belongs to $\cup_{t=0}^r \Lambda_t$. Its multiplicity as an eigenvalue of $A$ is given by
\[
m(\l) ~=~ \sum_{t: \l \in \Lambda_t} \({n \choose t} - {n \choose {t-1}}\).
\]
Here we define ${n \choose {-1}}$ to be $0$.

\item Let $f \in W_{t}(\l)$. Then $f$ vanishes on $B(n,t-1) = \cup_{i=0}^{t-1} S(n,i)$. Furthermore, the restriction of $f$ to each of the Hamming spheres $S(n,i)$, $t \le i \le r$ lies in the $t^{th}$ eigenspace, $V_t$ of $S(n,i)$. In particular, the restrictions of functions in $W_{t}(\l)$ to $S(n,t)$ span the $t^{th}$ eigenspace $V_t$ of $S(n,t)$. Moreover, a function $f \in W_{t}(\l)$ is determined by its restriction to $S(n,t)$.

\end{enumerate}

\ethm

\noi Some comments. 

\myblt The first two claims of the theorem describe the eigenvalues of $A$, and the third claim describes its eigenfunctions. In particular, the third claim of the theorem shows that there is a natural isomorphism beween the eigenspace $W_{t}(\l)$ and the $t^{th}$ eigenspace $V_t$ of $S(n,t)$.

\myblt It will be observed in the proof of the theorem, that if an eigenvalue $\l$ of $A$ belongs to two sets $\Lambda_{t_1}$ and $\Lambda_{t_2}$, with $t_1 \not = t_2$, then the eigenspaces $W_{t_1}(\l)$ and $W_{t_2}(\l)$ are orthogonal. Hence the third claim of the theorem implies that if $\l$ belongs to precisely $k$ sets $\Lambda_{t_i}$ for some $t_1 < t_2 < \ldots < t_k$, then there is a natural isomorphism between the space of all eigenfunctions of $A$ corresponding to the eigenvalue $\l$, and the direct sum $\oplus_{i=1}^k V_{t_i}$, where $V_{t_i}$ is the $t_i^{th}$ eigenspace of $S\(n,t_i\)$.

\myblt A Hamming ball is a subset of the Hamming cube which shares some properties of Hamming spheres. This is demonstrated in its spectral properties, in that its eigenvalues are closely related to the roots of the zonal spherical functions of Hamming cubes of varying dimensions, and the restrictions of its eigenfunctions to its constituent spheres lie in appropriate eigenspaces of these spheres.

\noi Coming back to the question of the maximal eigenvalue of a Hamming ball of radius $r$, we have the following claim.

\cor
\label{cor:max eigenvalue}
The maximal eigenvalue of $A$ is $\l = n-2x$, where $x$ is the first root of the Krawtchouk polynomial $K^{(n)}_{r+1}$. Its multiplicity is $1$, and the corresponding eigenfunction is a positive spherical function on $B$.
\ecor

\subsubsection{A more general setting - union of concentric Hamming spheres}
Theorem~\ref{thm:main} can be extended to cover a more general setting. Let $0 \le r_1 < r_2 \le n/2$. Let $B = B\(n,r_1,r_2\) = \cup_{i=r_1}^{r_2} S(n,i)$. Let $A = A\(n,r_1,r_2\)$ be the adjacency matrix of the subgraph of $\H$ induced by the vertices in $B$. We will describe the eigenvalues and the eigenfunctions of $A$ in the following result, which we present without proof, since the proof of Theorem~\ref{thm:main} extends essentially verbatim to this setting. 

\noi Let us first define a family of matrices and consider the eigenvalues of these matrices. For $0 \le t \le r_2$, we write $\t$ for $\max\{t,r_1\}$.

\dfn
\label{dfn:Matrices M}
For $0 \le t \le r_2$, we define the matrix $M_t = M_{n,r_1,r_2,t}$ to be the $\(r_2 - \t + 1\) \times \(r_2 - \t + 1\)$ matrix given by
\[
M_t ~=~ \(\begin{array}{ccccc} 0 & \sqrt{\l_{\t-t+2}} & 0 & \ldots & 0 \\ \sqrt{\l_{\t-t+2}} & 0 & \sqrt{\l_{\t-t+3}} & 0 & \ldots \\ \vdots & \vdots & \vdots & \vdots & \vdots \\ 0 & \ldots & \sqrt{\l_{r_2-t}} & 0 & \sqrt{\l_{r_2-t+1}} \\ 0 & \ldots & 0 & \sqrt{\l_{r_2-t+1}} & 0 \end{array} \),
\]
with $\l_k = (k-1)(n-2t-k+2)$.

\noi The matrix $M_t$ is symmetric and tridiagonal. Moreover, its off-diagonal elements are strictly positive. Hence $M_t$ has $r_2 - \t + 1$ simple real eigenvalues (see e.g. \cite{Chihara}, Theorem 4.4). We denote this set of eigenvalues by $\Lambda_t$.
\edfn

\noi Note that for all $0 < t \le r_1$, the matrices $M_t$ are $\(r_2 - r_1 + 1\) \times \(r_2 - r_1 + 1\)$ matrices. However, the entries in these matrices depend on $t$, and hence the sets $\Lambda_t$ are in general different (see the proof of Corollary~\ref{cor:incidence}). 

\noi We have the following claim.

\thm
\label{thm:general}

\begin{enumerate}

\item For each choice of $t = 0,...,r_2$ and for each of the $r_2-\t+1$ distinct real numbers $\l \in \Lambda_t$, the matrix $A$ has an eigenspace $W_{t}(\l)$ of dimension ${n \choose t} - {n \choose {t-1}}$ corresponding to the eigenvalue $\l$.

\item Let $f \in W_{t}(\l)$. Then the restriction of $f$ to each of the Hamming spheres $S(n,i)$, $t \le i \le r$ lies in the $t^{th}$ eigenspace $V_t$ of $S(n,i)$. Furthermore, $f$ vanishes on $\cup_{i=0}^{\t-1} S(n,i)$. The restrictions of functions in $W_{t}(\l)$ to $S\(n,\t\)$ span the $t^{th}$ eigenspace $V_t$ of $S\(n,\t\)$. Moreover, a function $f \in W_{t}(\l)$ is determined by its restriction to $S\(n,\t\)$.

\item A real number $\l$ is an eigenvalue of $A$ if and only if it belongs to $\cup_{t=0}^{r_2} \Lambda_t$. Its multiplicity as an eigenvalue of $A$ is given by
\[
m(\l) ~=~ \sum_{t: \l \in \Lambda_t} \({n \choose t} - {n \choose {t-1}}\).
\]
Here we define ${n \choose {-1}}$ to be $0$.

\end{enumerate}

\ethm

\noi We remark that for $r_1 = 0$ and $r_2 = r$ we recover Theorem~\ref{thm:main}, since in this case the eigenvalues of the matrix $M_t$ are given by the set $2 \cdot {\cal R}(n-2t,r-t+1) - (n-2t)$, see the proof of Lemma~\ref{lem:simple eigenvalues} below.

\subsection{Applications}

\subsubsection{Incidence matrices}

\noi Let us consider a special case of Theorem~\ref{thm:general}, in which $r_1 = r-1$ and $r_2 = r$, for some $1 \le r \le n/2$. That is, we deal with union of two adjacent concentric Hamming spheres. Note that the matrix $A$ has the following form 
\[
A ~=~ \(\begin{array}{ccccc} 0 & C & \\ C^t & 0 \end{array} \),
\]
where $C$ is the ${n \choose {r-1}} \times {n \choose r}$ {\it incidence matrix}, in which the rows are indexed by sets of size $r-1$ and the columns by sets of size $r$, and two sets are connected if the smaller one is contained in the larger one. Since (\cite{Gottlieb}) this matrix is known to be full-rank, we already know that the rank of $A$ is twice the rank of $C$, that is $2 \cdot {n \choose {r-1}}$.

\noi Theorem~\ref{thm:general} allows to obtain more information about $A$. In fact, the following result describes all the eigenvalues and the eigenvectors of $A$.

\cor 
\label{cor:incidence}
\begin{itemize}

\item For each $0 \le t \le r-1$, the matrix $A$ has a pair of eigenvalues $\pm \gamma_t$, with $\gamma_t = \sqrt{(r-t)(n-r-t+1)}$. The multiplicity of these eigenvalues is ${n \choose t} - {n \choose {t-1}}$. In addition $A$ has eigenvalue $0$ with multiplicity ${n \choose r} - {n \choose {r-1}}$.

\item Given $0 \le t \le r-1$, and $\l \in \{-\gamma_t, \gamma_t\}$, the eigenspace $W_{t}(\l)$ of $A$ which corresponds to the eigenvalue $\l$ can be described as follows. Let $f \in W_{t}(\l)$. Let $f_{r-1}$ and $f_r$ be the restrictions of $f$ to $S(n,r-1)$ and $S(n,r)$ respectively. Then $f_{r-1}$ can be chosen arbitrarily as a function in the $t^{th}$ eigenspace $V_t$ of $S\(n,r-1\)$, in which case the function $f_r$ is determined by the equation $f_r(x) = \lambda \cdot \sum_{|y| = r-1, y \subset x} f_{r-1}(y)$. Furthermore, $f_r$ lies in the $t^{th}$ eigenspace $V_t$ of $S\(n,r\)$.

\item The eigenspace $W(0)$ which corresponds to the eigenvalue $0$ contains the functions which vanish on $S(n,r-1)$ and whose restriction to $S(n,r)$ lies in the $r^{th}$ eigenspace $V_r$ of $S\(n,r\)$.
\end{itemize}
\ecor

\prf The claim of follows immediately from Theorem~\ref{thm:general}, since for $0 \le t \le r-1$, the matrices $M_t$ in this case are $2 \times 2$ matrices, $M_t ~=~ \(\begin{array}{ccccc} 0 & \gamma_t & \\ \gamma_t & 0 \end{array} \)$. Furthermore, the matrix $M_r$ is the $1 \times 1$ zero matrix.
\eprf 

\xam 
\label{xam:r=1}
Let $r=1$. Then the matrix $A$ is the following $(n+1) \times (n+1)$ matrix. 
\[
A ~=~ \(\begin{array}{ccccc} 0 & 1 & \cdots & 1 \\ 1 & 0 & \cdots & 0 \\ \vdots & \vdots & \cdots & \vdots \\ 1 & 0 & \cdots &  0  \end{array} \).
\]
Corollary~\ref{cor:incidence} gives that in this case $A$ has eigenvalues $\pm \sqrt{n}$ with mulitplicity $1$ each, and eigenvalue $0$ with multiplicty $n-1$. This can be easily confirmed by direct computation. 
\exam

\rem 
\label{rem:adjacent spheres - easy}

\noi Let us sketch out an alternative way to obtain the results of Corollary~\ref{cor:incidence}. It is easy to see that $A^2$ is a block-diagonal matrix with two blocks. One of these blocks is an ${n \choose {r-1}} \times {n \choose {r-1}}$ matrix, which is a linear combination of the first two adjacency matrices of $S(n,r-1)$. Similarly, the second block is an ${n \choose {r}} \times {n \choose {r}}$ matrix which is a linear combination of the first two adjacency matrices of $S(n,r)$. This allows to find the eigenvalues and the eigenfunctions of $A^2$ using the theory of association schemes (specifically of the Johnson association scheme related to the Hamming sphere). The eigenvalues and the eigenfunctions of $A$ are now easy to determine, using the fact that $A$ is a matrix of a bipartite graph.
\erem

\subsubsection{Maximal eigenvalues and minimal fractional edge boundaries of large sets}

\noi We will use Theorem~\ref{thm:main} to show that Hamming balls have essentially the largest maximal eigenvalue among all subsets of the Hamming cube with the same cardinality, even when the sets in question are large, with cardinality which is an arbitrary subconstant fraction of the whole cube. 

\noi We recall relevant notions and introduce notation. For a subset $A \subseteq \H$, let $\l(A)$ denote the maximal eigenvalue of the subgraph of the Hamming cube induced by the vertices in $A$. Define the fractional edge boundary size (\cite{sam-log-sob})  of $A$ by
\[
|\partial^{\ast} A| ~=~ \min \left\{\frac{1}{2^{n+1}}  \sum_{x,y:\,|x-y|=1} \big(f(x)-f(y)\big)^2 ~\bigg |~ f:\H \rarrow \R,~\mathrm{supp}(f) \subseteq A,~ \frac{1}{2^n} \sum_x f^2(x) = 1\right\}.
\]

\noi It is easy to see that the $|\partial^{\ast}(A)| = n - \l(A)$. 

\noi The question of the largest maximal eigenvalue, equivalently the minimal fractional edge boundary size, for subsets of $\H$ with given cardinality was raised in \cite{ft}. A partial answer to this question was given in \cite{sam-log-sob}, where the following lower bound on the fractional edge boundary size was given. Let $H(x) = x \log_2\(\frac{1}{x}\) + (1-x) \log_2\(\frac{1}{1-x}\)$ be the binary entropy function. This function (defining it to be zero at $0$) is strictly increasing on $\left[0,1/2\right]$, and it takes $\left[0,1/2\right]$ to $\left[0,1\right]$. Let A be a subset of $\H$ of cardinality $1 \le s \le 2^{n-1}$. Then 
\beqn
\label{ineq:modLS}
|\partial^{\ast} A| ~\ge~ n \cdot \(1 - 2\sqrt{H^{-1}\(\frac{\log_2(s)}{n}\) \(1 - H^{-1}\(\frac{\log_2(s)}{n}\)\)}\).
\eeqn

\noi To connect notions, it seems useful to compare this inequality with the 'usual' edge-isoperimetric inequality \cite{Harper, Hart} in the Hamming cube. Recall that the edge boundary size $|\partial A|$ of a set $A$ is the number of edges between $A$ and $A^c$. The edge-isoperimetric inequality then states that $\frac{|\partial A|}{|A|} \ge  \log_2\(\frac{2^n}{|A|}\) = n \cdot \(1 - \frac{\log_2(s)}{n}\)$, where $s := |A|$. This is tight if $A$ is a subcube of the Hamming cube. Observe that in order to determine the fractional edge boundary size of a set, we compute the minimum of a certain notion of variation over a family of functions supported in this set. Since this family contains the (normalized) characteristic function of $A$, whose variation is proportional to the edge boundary size of $A$, it is easy to verify that $|\partial^{\ast} A| \le \frac{|\partial A|}{|A|}$. 
Furthermore, it is easy to see that $|\partial^{\ast} A| = \frac{|\partial A|}{|A|} = \log_2\(\frac{2^n}{|A|}\)$ if $A$ is a subcube. 
Hence (\ref{ineq:modLS}) provides a lower bound on $\frac{|\partial A|}{|A|}$, which is in general weaker than the edge-isoperimetric inequality (in fact, it is easy to see that $2\sqrt{H^{-1}(x) \(1 - H^{-1}(x)\)} > x$, for any $0 < x < 1$). In particular, (\ref{ineq:modLS}) is not tight for subcubes.
 
\noi For $1 \le s \le 2^{n-1}$, let $\Lambda(s)$ denote the largest maximal eigenvalue among all subsets of the Hamming cube of cardinality $s$. Let $\Delta(s)$ be the smallest minimal fractional edge boundary size among these subsets. Clearly, $|\Delta(s)| = n - \Lambda(s)$. Introducing the parametrization $r = r(s) = n H^{-1}\(\frac{\log_2(s)}{n}\)$, the preceding discussion shows that, at least when $s$ is a power of $2$, 
\[
n H\(\frac{r}{n}\) ~\le~ \Lambda(s) ~\le~ 2\sqrt{r(n-r)}.
\] 
The lower bound is given by subcubes of cardinality $s = 2^{n H\(\frac{r}{n}\)}$, and the upper bound by (\ref{ineq:modLS}).

\noi Tighter lower bounds on $\Lambda(s)$ were obtained in \cite{ft, BLL} by analyzing the maximal eigenvalues of Hamming balls which, for many values of $s$, turn out to be larger that those of subcubes of the same (or similar) size. In particular, the results in \cite{BLL} imply 
\[
\Lambda(s) ~\ge~ 2\sqrt{r \(n+1-r\)} \cdot \(1 - O\(\sqrt{\frac{\log\(r\)}{r}}\)\),
\] 
when $r \rarrow \infty$ (similar estimates were obtained in \cite{ft}). 

\noi This implies that for $r \rarrow \infty$ the inequality $\Lambda(s) \le 2\sqrt{r(n-r)}$ is essentially tight, that is $\(1-o_r(1)\) \cdot 2\sqrt{r(n-r)} \le \Lambda(s) \le 2\sqrt{r(n-r)}$, and that Hamming balls of growing radius have essentially the largest maximal eigenvalue among all subsets of the Hamming cube with the same cardinality. For a fixed radius $r$, the situation is less satisfactory, in that the 'error term' is of the same magnitude as the 'main term' $2\sqrt{r(n-r)}$. To complement this result, it is shown in \cite{BLL} that for every fixed $r$ there is an $n_0 = n_0(r)$ so that for any $n \ge n_0$, the Hamming ball of radius $r$ has {\it the largest} maximal eigenvalue among all subsets of the $n$-dimensional Hamming cube with the same cardinality.

\noi Passing to the fractional edge boundary size, we have, for $1 \le s \le 2^{n-1}$ and $r = n H^{-1}\(\frac{\log_2(s)}{n}\)$, that 
\[
n - 2\sqrt{r(n-r)} ~\le~ \Delta(s) ~\le~ n - 2\sqrt{r \(n+1-r\)} \cdot \(1 - O\(\sqrt{\frac{\log\(r\)}{r}}\)\). 
\]

\noi This shows that $n - 2\sqrt{r(n-r)} \le \Delta(s) \le \(1+o_n(1)\) \cdot \(n - 2\sqrt{r(n-r)}\)$ for $r \le \frac n2 - \widetilde{\Omega}\(n^{3/4}\)$, and hence Hamming balls of radius $r$ have essentially the smallest fractional edge boundary size.

\noi For larger sets the situation is less clear. Consider the case $s = 2^{n-1}$. It is not hard to see\footnote{This is a familiar fact, which can be proved by a simple application of Fourier analysis on the boolean cube.} that among all sets of cardinality $s$, the set $A$ with the smallest fractional edge boundary is the $(n-1)$-dimensional subcube (rather than a Hamming ball of radius $\frac{n-1}{2}$, assuming $n$ is odd). In this case estimating $\Delta(s)$ by $n - 2\sqrt{r(n-r)}$ errs by a constant multiplicative factor. In fact, we have $\Delta\(s\) = |\partial^{\ast} A| = 1$, while $r = r(s) \approx n/2 - \sqrt{\(\ln(2)/2\) \cdot n}$ (to see this, recall that $H\(1/2 - x\) = 1 - 2/\ln(2) \cdot x^2 + o\(x^2\)$), and hence $n - 2\sqrt{r(n-r)} \approx \ln(2)$.    

\noi More generally, let $s = 2^{n-o(n)}$. Then $r(s) \approx n/2 - \sqrt{\(\ln\(2^n/s\)/2\) \cdot n}$, giving the lower bound $\Delta(s) ~\ge~ n - 2\sqrt{r(n-r)} ~\approx~ \ln\(2^n/s\)$. On the other hand, the example of Hamming subcubes implies the upper bound $\Delta(s) \le \log_2\(2^n/s\)$, at least when $s$ is a power of two. There is a multiplicative gap of $1/\ln(2)$ between these two bounds and, as observed above, the upper bound is tight for $s = 2^{n-1}$. Extending this observation, it was conjectured in \cite{sam:faber krahn large sets} that $\Delta(s) \ge \(1-o_n(1)\) \cdot \log_2\(2^n/s\)$ for $2^n/n \le s \le 2^{n-1}$.

\rem
\label{rem:open q}
More precisely and more generally, the results in \cite{sam-log-sob}, as well as the logarithmic Sobolev inequality \cite{Gross}, imply that $\Delta(s) \ge \ln\(2^n/s\)$ for all $1 \le s \le 2^{n-1}$. Equivalently, $\Lambda(s) \le n - \ln\(2^n/s\)$. This provides a positive answer to Question~27 in \cite{BLL}.
\erem

\subsubsection*{Our results}

\noi Corollary~\ref{cor:max eigenvalue} in conjunction with known properties of Krawtchouk polynomials and their roots can be applied to tighten the bounds on the maximal eigenvalues of Hamming balls, and consequently on the functions $\Lambda$ and $\Delta$. We show the following two results.

\cor 
\label{cor:max-ev}

\begin{itemize}

\item Let $1 \le s \le 2^{n-1}$ and let $r = n H^{-1}\(\frac{\log_2(s)}{n}\)$. Then
\[
\Lambda(s) ~\ge~ 2\sqrt{r \(n-r\)} - O_{r \rarrow \infty}\(r^{-1/6} \cdot \sqrt{n}\).
\]

\item Let $s = 2^{n-o(n)}$. Then 
\[
\Delta(s) ~\le~ \ln\(\frac{2^n}{s}\) \cdot \(1 + o_{s/2^n \rarrow 0}(1)\).
\]

\end{itemize}
\ecor

\noi Some comments. 

\myblt The error term in the first claim of the corollary is smaller than in the above mentioned bound of \cite{BLL}: $O\(r^{-1/6} \cdot \sqrt{n}\)$ rather than $O\(\sqrt{\log(r)} \cdot \sqrt{n}\)$.

\myblt The second claim of the corollary disproves the above mentioned conjecture of \cite{sam:faber krahn large sets}.

\myblt The second claim, together with the discussion above, implies that Hamming balls of cardinality $s$ have essentially the largest maximal eigenvalue (respectively the smallest fractional boundary size) among all subsets of the Hamming cube with the same cardinality, when $s$ is an arbitrary subconstant fraction of $2^n$.

\section{Proofs}

\subsection{Proof of Proposition~\ref{pro:semi-symmetric}}

\noi Let $0 \le j \le i \le n/2$. Let $V_j$ be the $j^{th}$ eigenspace of $S(n,i)$. Let $y \in \H$, $|y| = j$. Recall that we denote by $S_{y,i}$ the space of all semi-symmetric functions around $y$ in $S(n,i)$. Therefore the space of zonal semi-symmetric functions around $y$ in $V_j$ is given by $V_{y,j} = V_j \cap S_{y,i}$. Our goal is to show that this space is $1$-dimensional.

\noi We start with studying the space $S_{y,i}$. Note this is a $(j+1)$-dimensional subspace of functions on $S(n,i)$, spanned by the characteristic functions of the subsets $B_k = B_{y,k} = \{x \in S(n,i),~|x \cap y| = k\}$, $k=0,...,j$. Let us construct an additional basis of $S_{y,i}$. For $0 \le k \le j$, let $G_k = G_{y,k} = \sum_{z \subseteq y, |z| = k} g_z$ (the functions $g_z$ are defined in Definition~\ref{dfn:functions g}). Observe that $G_k = \sum_{\ell=0}^j {{\ell} \choose k} \cdot 1_{B_{\ell}}$, where we define the binomial coefficient ${b \choose a}$ to be $0$ if $b < a$. Hence, the functions $\{G_k\}_{k=0}^j$ are in $S_{y,i}$, and moreover it is easy to see that they are linearly independent. It follows that they form another basis of $S_{y,i}$.

\noi It useful to note that for any $z \subseteq y$, the function $G_{|z|}$ is proportional to the orthogonal projection of $g_z$ on $S_{y,i}$. In fact, let $z_1 \not = z_2$ be two subsets of $y$ with the same cardinality $k$. Let $f \in S_{y,i}$. We claim that $g_{z_1} - g_{z_2} \perp f$. To see that, let $\sigma \in S_n$ take $z_1$ to $z_2$ while preserving $y$. Let $T_{\sigma}$ (see Section~\ref{subsubsec:doubly transitive}) be the corresponding linear operator on the functions on $S(n,i)$. Note that $T_{\sigma}:~g_{z_1} \rarrow g_{z_2}$, while $T_{\sigma} f = f$. Note also that $T_{\sigma}$ preserves the inner product in $S(n,i)$. Hence $\<g_{z_1}, f\> = \<T_{\sigma} g_{z_1}, T_{\sigma} f\> = \<g_{z_2},f\>$. It follows that for $z \subseteq y$, the orthogonal projection $P_{y,i}$ of $g_z$ on $S_{y,i}$ depends only on $|z|$, and hence, writing $k = |z|$, we get $P_{y,i} \, g_z = \frac{1}{{j \choose k}} \cdot \(P_{y,i} ~\sum_{w \subseteq y, |w| = k} g_w\) = \frac{1}{{j \choose k}} \cdot P_{y,i}\(G_k\) = \frac{1}{{j \choose k}} \cdot G_k$.

\noi Let us define $U_{y,j} = S_{y,i} \cap \<\{g_z\}_{z \subsetneq y}\>^{\perp}$. By the discussion above, this is equivalent to $U_{y,i} = S_{y,i} \cap \<G_0,...,G_{j-1}\>^{\perp}$. This means that the space $U_{y,i}$ is $1$-dimensional. We will show that $U_{y,j} = V_{y,j}$ and this will complete the proof of the proposition.

\noi We first claim that $U_{y,j} = S_{y,i} \cap \<\{g_z\}_{|z| < j}\>^{\perp}$. This will follow if we show that for each $z \in \H$, $|z| < j$, holds $P_{y,i}\(g_z\) \in \<G_0,...,G_{j-1}\>$. By using symmetry as above, it is easy to see that $P_{y,i}\(g_z\)$ depends only on $|z|$ and on $|z \cap y|$. A simple consequence of this is that  $P_{y,i}\(g_z\)$ is proportional to $\sum_{\ell=0}^j {{\ell} \choose {|z \cap y|}}{{i - \ell} \choose {|z| - |z \cap y|}}  \cdot 1_{B_{\ell}}$, where we again define ${b \choose a}$ to be $0$ if $b < a$. Observe that the vector of coefficients $F(\ell) = {{\ell} \choose {|z \cap y|}}{{i - \ell} \choose {|z| - |z \cap y|}}$ is given by a restriction to the integer points $0,1,...,j$ of a degree-$|z|$ real valued polynomial in $\ell$, and that the vectors of coefficients of the functions $G_k = \sum_{\ell=0}^j {{\ell} \choose k} \cdot 1_{B_{\ell}}$, $k=0,...,j-1$ span all polynomials of degree up to $j-1$ restricted to $0,1,...,j$. Hence $P_{y,i}\(g_z\) \in \<G_0,...,G_{j-1}\>$.

\noi Next, observe that $S_{y,i} = \<G_0,...,G_j\> \subseteq \<\{g_z\}_{|z| \le j}\>$. It follows (see the description of the eigenspaces of $S(n,i)$ in Section~\ref{subsubsec:doubly transitive}) that
\[
U_{y,j} ~=~ S_{y,i} \cap \<\{g_z\}_{|z| < j}\>^{\perp} ~=~ \Big(S_{y,i} \cap \<\{g_z\}_{|z| \le j}\>\Big) \cap \<\{g_z\}_{|z| < j}\>^{\perp} ~=~
\]
\[
S_{y,i} \cap \Big(\<\{g_z\}_{|z| \le j}\> \cap \<\{g_z\}_{|z| < j}\>^{\perp}\Big) ~=~ S_{y,i} \cap \Big(U_j \cap U_{j-1}^{\perp}\Big) ~=~ S_{y,i} \cap V_j ~=~ V_{y,j}.
\]

\noi To complete the proof of the proposition, it remains to show that $V_j = \<V_{y,j}\>_{|y| = j}$. Let $E_j$ be the orthogonal projection on $V_j$ (see Section~\ref{subsubsec:doubly transitive}). We claim that the function $f = E_j 1_{B_{y,j}}$ is in $V_{y,j}$. Since, by definition of $E_j$, $f$ is in $V_j$, we need to verify that $f \in S_{y,i}$. By Lemma~\ref{lem:y-invariant}, it suffices to show that $f$ is preserved by isometries which preserve $y$. Let $\sigma \in S_n$ be a $y$-preserving permutation. Then, since $1_{B_{y,j}}$ is in $S_{y,i}$ and hence is preserved by $T_{\sigma}$, and since $E_j$ commutes with $T_{\sigma}$, we have
\[
T_{\sigma} f ~=~ T_{\sigma} E_j 1_{B_{y,j}} ~=~ E_j  T_{\sigma} 1_{B_{y,j}} ~=~ E_j 1_{B_{y,j}} ~=~ f.
\]

\noi Next, observe that $1_{B_{y,j}} =  \{x \in S(n,i),~|x \cap y| = j\}$ is precisely the function $g_y$, as defined in Definition~\ref{dfn:functions g}. Recall also (see Section~\ref{subsubsec:doubly transitive}) that the functions $\{g_y\}_{|y| = j}$ span a subspace $U_j$, which contains $V_j$. Hence
\[
V_j ~=~ E_j U_j ~=~ E_j \<\{g_y\}\>_{|y| = j}  ~=~ E_j \<1_{B_{y,j}}\>_{|y| = j} ~\subseteq~ \<V_{y,j}\>_{|y| = j}.
\]

\noi This completes the proof of the proposition. \eprf

\noi For future use, let us record a useful fact which stems from the argument above. Let $f \in V_{y,j}$ be a non-zero function. We claim that $f$ does not vanish on $B_j = B_{y,j}$. To see that, note that $1_{B_j} = G_j$, and hence the functions $G_0,...,G_{j-1}$ together with $1_{B_j}$ span $S_{y,i}$. Furthermore, we showed that $V_{y,j} = U_{y,i} = S_{y,i} \cap \<G_0,...,G_{j-1}\>^{\perp}$. It follows that if $f$ vanishes on $B_j$, it is orthogonal to the whole space $S_{y,i}$, and hence has to be the zero function.

\subsection{Proof of Theorem~\ref{thm:main}}

\noi Recall that $B = B(n,r)$ is the Hamming ball of radius $r$ around $0$ in $\H$, and that $A = A(n,r)$ is the adjacency matrix of the subgraph of $\H$ induced by the vertices in $B$. We start with a description of some relevant properties of the space of functions on $B$. The main observation is that since $B$ is invariant under a large group of isometries of $\H$, namely the permutation group $S_n$, it inherits many of the convenient properties of the Hamming sphere.

\noi For $y \in B$, we define $S^{(B)}_y$ to be the vector space of all semi-symmetric functions around $y$ on $B$. These are the functions on $B$ whose value at a point $x \in B$ depends only on the Hamming distances $|x| = |x-0|$ and $|x-y|$. Recall that we denote by $S_y$ the vector space of all semi-symmetric functions around $y$ on $\H$, and note that $S^{(B)}_y$ is obtained by restricting the functions in $S_y$ to $B$.

\noi The following property of $S^{(B)}_y$ will play an important role in the argument below.

\lem
\label{lem:V-invariant}
The space $S^{(B)}_y$ is an invariant subspace of $A$.
\elem
\prf
We first consider functions on $\H$ and prove that $S_y$ is an invariant subspace of $A_1$, the adjacency matrix of $\H$. Let $f \in S_y$. We claim that $Af \in S_y$. Let $x_1$ and $x_2$ be two points in $\H$ with the same Hamming distance from $0$ and $y$. Then there is a permutation $\sigma$ of $\H$ which preserves $y$ and which  takes $x_1$ to $x_2$. Recall (see Section~\ref{subsubsec:doubly transitive}) that $A_1$ commutes with the isometries $\{T_{\alpha}\}$ of $\H$, and that for any permutation $\alpha$ preserving $y$ we have $T_{\alpha} f = f$. Hence, writing $\delta_x$ for the characteristic function of $x$,
\[
(A_1 f)\(x_2\) ~=~ \<A_1 f, \delta_{x_2}\> ~=~ \<A_1 f, T_{\sigma} \delta_{x_1}\> ~=~ \<T_{\sigma^{-1}} A_1 f, \delta_{x_1}\> ~=~ \<A_1 T_{\sigma^{-1}} f, \delta_{x_1}\> ~=~ \<A_1 f, \delta_{x_1}\> ~=~ f\(x_1\).
\]

\noi Next, we consider functions on $B$. Let $f \in S^{(B)}_y$. Let $g$ be the extension of $f$ to the whole $\H$ obtained by defining $g = f$ on $B$ and $0$ outside $B$. Note that $Af$ equals to the restriction of $A_1 g$ to $B$. Since $f \in S^{(B)}_y$, it is easy to see that $g$ is in $S_y$. By the preceding argument, also $A_1 g$ is in $S_y$, and hence $Af$ is in $S^{(B)}_y$, as the restriction of $A_1 g$ to $B$.
\eprf

\noi We proceed with the proof of the theorem. For $y \in B$, we define an $A$-invariant subspace of functions which depends on $y$, and study the restriction of $A$ to this subspace.

\dfn
\label{dfn: matrix A}
Let $y \in B$. We define a subspace $W_y$ of functions on $B$ as follows:
\[
W_y ~=~ S^{(B)}_y \cap \<\left\{S^{(B)}_z\right\}_{|z| < |y|}\>^{\perp}.
\]
\edfn

\noi Here the span of an empty family of spaces is defined to be $\{0\}$. Hence if $y = 0$, we have $W_y = S^{(B)}_y$. For any $y \in B$, $W_y$ is an invariant subspace of $A$, as an intersection of invariant subspaces. 

\noi We have the following useful lemma. Recall that we denote by $S_{y,i}$ the space of all semi-symmetric functions around $y$ in $S(n,i)$. For $|y| = t$ we define $V_{y,i}$ as the space of zonal semi-symmetric functions around $y$ in the $t^{th}$ eigenspace $V_t$ of $S(n,i)$.

\lem
\label{lem:key}
Let $|y| = t$. For $0 \le i \le r$, let $W_{y,i}$ be the restriction of $W_y$ to $S(n,i)$. Then $W_{y,i} = \{0\}$ for $i < t$ and $W_{y,i} = V_{y,i}$, for $t \le i \le r$. \elem

\prf
Let $0 \le i \le r$. We claim that $W_{y,i}  = S_{y,i} \cap \<\{S_{z,i}\}_{|z| < t}\>^{\perp}$. In fact, let $f \in W_y$, and let $f_i$ be the restriction of $f$ to $S(n,i)$. Note that  $f \in S^{(B)}_y$ implies $f_i \in S_{y,i}$. We claim that $f_i$ is orthogonal to $\<\{S_{z,i}\}_{|z| < |t|}\>$. In fact, for any $z$ the restriction of $S^{(B)}_z$ to $S(n,i)$ is $S_{z,i}$. Furthermore, if we view $S_{z,i}$ as a space of functions on $B$ which vanish outside $S(n,i)$, we have $S_{z,i} \subseteq S^{(B)}_z$. Hence $f \perp S^{(B)}_z$ implies $f \perp S_{z,i}$, which in its turn implies $f_i \perp S_{z,i}$. It follows that $W_{y,i} \subseteq S_{y,i} \cap \<\{S_{z,i}\}_{|z| < t}\>^{\perp}$. Conversely, let $f_i \in S_{y,i} \cap \<\{S_{z,i}\}_{|z| < t}\>^{\perp}$. Viewing $f_i$ as a function on $B$ which vanishes outside $S(n,i)$, note that $f_i \in S^{(B)}_y$, and moreover, that it is orthogonal to $\<\left\{S^{(B)}_z\right\}_{|z| < t}\>$. Hence $f_i \in W_y$ and therefore, viewed as a function on $S(n,i)$, it is also in $W_{y,i}$, and we have an inclusion in the opposite direction as well.

\noi Next, we claim that for any $i \ge t$ we have $\<\{S_{z,i}\}_{|z| < t}\> = \<\{g_z\}_{|z| < t}\>$ (where the functions $g_z$ are defined in Definition~\ref{dfn:functions g}). To see that, first note that for any $z \in \H$, $|z| \le i$, we have $g_z \in S_{z,i}$, and hence the first space contains the second one. On the other hand, as shown in the proof of Proposition~\ref{pro:semi-symmetric}, for any $z \in \H$, $|z| \le i$, holds $S_{z,i} \subseteq \<\{g_w\}_{|w| \le |z|}\>$. The inclusion in the second direction follows.

\noi We can now complete the proof of the lemma. First, let $0 \le i < t$. Then
\[
W_{y,i} ~=~ S_{y,i} \cap \<\{S_{z,i}\}_{|z| < t}\>^{\perp} ~\subseteq~ S_{y,i} \cap \<\{S_{z,i}\}_{|z| \le i}\>^{\perp} ~=~ S_{y,i} \cap \<\{g_z\}_{|z| \le i}\>^{\perp} ~=~ \{0\}.
\]
For the last step observe that the functions $\<\{g_z\}_{|z| = i}\>$ span all functions on $S(n,i)$. 

\noi Let now $t \le i \le r$. In this case,
\[
W_{y,i} ~=~ S_{y,i} \cap \<\{S_{z,i}\}_{|z| < t}\>^{\perp} ~=~ S_{y,i} \cap \<\{g_z\}_{|z| < t}\>^{\perp} ~=~ V_{y,i},
\]
where the last equality is shown in the proof of Proposition~\ref{pro:semi-symmetric}.

\noi For future use, let us record a useful fact which follows from the argument above. Let $t \le i \le r$. Then the space $W_{y,i} = V_{y,i}$ considered as a space of functions on $B$, is a subspace of $W_y$.

\eprf

\cor
\label{cor:A restricted to V_y}

Let $|y| = t$. Then the dimension of $W_y$ is $r-t+1$. Furthermore, $W_y$ has a basis consisting of zonal semi-symmetric functions $f_{y,t},...,f_{y,r}$ with $f_{y,i} \in V_{y,i}$ such that $f_{y,i}(x) = 1$ for $x \in S(n,i)$ with $y \subseteq x$ (and this determines $f_{y,i}$ uniquely). The restriction $A_y$ of $A$ to $W_y$ has the following form, in terms of the basis $f_{y,t},...,f_{y,r}$:
\[
A_y ~=~ \(\begin{array}{ccccc} 0 & \beta_{t+1} & 0 & \ldots & 0 \\ \gamma_t & 0 & \beta_{t+2} & 0 & \ldots \\ \vdots & \vdots & \vdots & \vdots & \vdots \\ 0 & \ldots & \gamma_{r-2} & 0 & \beta_r \\ 0 & \ldots & 0 & \gamma_{r-1} & 0 \end{array} \)
\]
Here $\beta_i = n - i + 1$ and $\gamma_i = \frac{(i-t+1)(n-t-i)}{n-i}$.
\ecor

\prf
By Lemma~\ref{lem:key}, $W_y = \oplus_{i=t}^r W_{y,i} = \oplus_{i=t}^r V_{y,i}$. By Proposition~\ref{pro:semi-symmetric}, the spaces $V_{y,i}$ are $1$-dimensional, and hence 
the dimension of $W_y$ is $r-t+1$. As observed in the comment following the proof of Proposition~\ref{pro:semi-symmetric}, a non-zero function in $V_{y,i}$ does not vanish on $x \in S(n,i)$ with $y \subseteq x$, and hence it is determined by normalising it to be $1$ on this set. Let $f_{y,i}$ be this function. Then $f_{y,t},...,f_{y,r}$ they form a basis of $W_y$.

\noi Let $A_y$ be the restriction of $A$ to $W_y$. Observe that in general $A$ acts on functions on $B$ in the following manner. Let $g:B \rarrow \R$. Then, for $x \in B$ we have $(Ag)(x) = \sum_{z: |x-z| = 1} g(y)$. The fact that $W_y$ is an invariant subspace of $A$ implies that there exist constants $\{\beta_i\}, \{\gamma_i\}$ so that 

\begin{itemize}

\item $A_y f_{y,t} = \gamma_t \cdot f_{y,t+1}$.

\item $A_y f_{y,r} = \beta_r \cdot f_{y,r-1}$.

\item For $t < i < r$, $A_y f_{y,i} = \beta_i \cdot f_{y,i-1} + \gamma_i \cdot f_{y,i+1}$.

\end{itemize}

\noi Note that the coefficients $\beta_i, \gamma_i$ should not depend on $y$ (for two points $y_1$ and $y_2$, consider an isometry taking $y_1$ to $y_2$).

\noi This already shows that the matrix $A_y$ has a tri-diagonal form, as claimed by the Corollary.  It remains to determine these coefficients. First, we claim that $\beta_i = n-i+1$, $i = t+1,...,r$. To see that, fix $i$ in this range, and consider the value of $A_y f_{y,i}$ at a point $x \in S(n,i-1)$ which contains $y$. On one hand $\(A_y f_{y,i}\)(x) = \sum_{z \in S(n,i), x \subseteq z} f_{y,i}(z) = n-i+1$. On the other hand, $\(A_y f_{y,i}\)(x) = \beta_i f_{y,i-1}(x) = \beta_i$.

\noi Next, we claim that $\gamma_i = \frac{(n-i-t)(i-t+1)}{n-i}$, for $i = t,...,r-1$. To show this, we first determine the value of $f_{y,i}(x)$ for $x \in S(n,i)$, $|x \cap y| = t-1$. As a zonal semi-spherical function, $f_{y,i}$ is orthogonal (see the proof of Proposition~\ref{pro:semi-symmetric}) to the function $G_{t-1} = 1_{B_{t-1}} + t \cdot 1_{B_t}$ in $S_{y,i}$. Here $B_j \subseteq S(n,i)$ is the set of all points $x$ for which $|x \cap y| = j$. Since $f_{y,i}$ is $1$ on $B_t$,  on $B_{t-1}$ it has to be equal to
\[
\alpha_{t-1} ~=~ -\frac{t |B_t|}{|B_{t-1}|} ~=~ -\frac{t \cdot {{n-t} \choose {i-t}}}{t \cdot {{n-t} \choose {i-t+1}}} ~=~ -\frac{i-t+1}{n-i}.
\]

\noi Fix $t \le i \le r-1$, and consider the value of $A_y f_{y,i}$ at a point $x \in S(n,i+1)$ which contains $y$. On one hand $\(A_y f_{y,i}\)(x) = \sum_{z \in S(n,i), z \subseteq x} f_{y,i}(z) = (i-t+1) + t \alpha_{t-1} = \frac{(n-i-t)(i-t+1)}{n-i}$. On the other hand, $\(A_y f_{y,i}\)(x) = \gamma_i f_{y,i+1}(x) = \gamma_i$.

\noi This concludes the proof of the Corollary.
\eprf

\noi The eigenvalues of the matrix $A_y$ are described in the following lemma. Recall that for $X \subseteq \R$ and $s, m \in \R$, we write $m X + s$ for the set $\{mx + s\}_{x \in X}$.

\lem
\label{lem:simple eigenvalues}
Let $A_y$ be the matrix defined in Corollary~\ref{cor:A restricted to V_y}. It has $r - t + 1$ simple eigenvalues given by the set $2 \cdot {\cal R}(n-2t,r-t+1) - (n-2t)$.
\elem

\prf

\noi For $k = t+1,...,r$, let $\l'_{k-t+1} = \beta_{k} \gamma_{k-1} = (k-t)(n-t-k+1)$. Write $m$ for $r - t + 1$, and let $D$ be an $m \times m$ diagonal matrix with $x_1,...,x_m$ on the diagonal, where $x_1 = 1$ and $x_k = \sqrt{\frac{\prod_{i=1}^{k-1} \gamma_{i+t-1}}{\prod_{i=2}^k \beta_{i+t-1}}}$, for $1 < k \le m$. It is easy to verify that the matrix $A'_y = D^{-1} A_y D$ is given by 
\[
A'_y ~=~ \(\begin{array}{ccccc} 0 & \sqrt{\l'_2} & 0 & \ldots & 0 \\ \sqrt{\l'_2} & 0 & \sqrt{\l'_3} & 0 & \ldots \\ \vdots & \vdots & \vdots & \vdots & \vdots \\ 0 & \ldots & \sqrt{\l'_{m-1}} & 0 & \sqrt{\l'_m} \\ 0 & \ldots & 0 & \sqrt{\l'_m} & 0 \end{array} \)
\]
In particular, the eigenvalues of $A_y$ are the same as those of $A'_y$.

\noi We use the following claim (\cite{Chihara}, Ex. 5.7). Let $P_0,...,P_m$ be a sequence of monic polynomials defined by the recurrence formula $P_k(x) = (x - c_k) P_{k-1}(x) - \l_k P_{k-2}(x)$, $k = 1,2,...$, with $P_{-1}(x) = 0$ and $P_0(x) = 1$. Here $\{c_k\}$ are real numbers, and $\{\l_k\}$ are strictly positive real numbers.
Let $M$ be the $m \times m$ matrix
\[
M ~=~ \(\begin{array}{ccccc} c_1 & \sqrt{\l_2} & 0 & \ldots & 0 \\ \sqrt{\l_2} & c_2 & \sqrt{\l_3} & 0 & \ldots \\ \vdots & \vdots & \vdots & \vdots & \vdots \\ 0 & \ldots & \sqrt{\l_{m-1}} & c_{m-1} & \sqrt{\l_m} \\ 0 & \ldots & 0 & \sqrt{\l_m} & c_m \end{array} \)
\]
Then the eigenvalues of $M$ are the zeroes of $P_m$.

\noi Consider the Krawtchouk polynomials $K_j = K^{(N)}_j = \sum_{\ell=0}^j (-1)^{\ell} {x \choose {\ell}} {{N-x} \choose {j-\ell}}$, $0 \le j \le N$ (see Section~\ref{subsubsec:doubly transitive}). These polynomials satisfy the recurrence formula $k K_k(x) = (N - 2x) K_{k-1}(x) - (N-k+2)  K_{k-2}(x)$ (see e.g., \cite{mrrw}). These polynomials are not monic. In fact, the main coefficient of $K_k$ is $\frac{(-2)^k}{k!}$. Let $P_k = \frac{k!}{(-2)^k} \cdot K_k$. The polynomials $\{P_k\}$ are monic and satisfy the recurrence $P_k(x) = \(x - \frac N2\) P_{k-1}(x) - \frac{(k-1)(N-k+2)}{4} P_{k-2}(x)$. Take $N = n-2t$, and $m = r - t + 1$. It follows that if in the matrix $M$ we take $c_1 = ... = c_m = \frac{n-2t}{2}$ and $\lambda_k = \frac{(k-1)(n-2t-k+2)}{4}$, for $k=2,...,m$, the obtained matrix has $m$ distinct eigenvalues, which are the roots of the Krawtchouk polynomial $K^{(n-2t)}_m$.

\noi It is easy to verify that the matrix $M$ we obtain specializing to these values of $\{c_k\}$ and $\l_k$ equals $\frac{n-2t}{2} \cdot I + \frac12 \cdot A'_y$, where $I$ is the $m \times m$ identity matrix. It follows that the eigenvalues of $A'_y$ are given by $2x_k - (n-2t)$, $k = 1,...,m$, where $\{x_k\}$ are the roots of $K^{(n-2t)}_m$. This completes the proof of the lemma.

\eprf

\noi Let us recap what we know so far. Let $0 \le t \le r$. Let $y \in S(n,t)$. Let $\l \in 2 \cdot {\cal R}(n-2t,r-t+1) - (n-2t)$. By Lemma~\ref{lem:simple eigenvalues}, $\l$ is a simple eigenvalue of $A_y$. Observe that $v_{\l} \in \R^{r-t+1}$ is a corresponding eigenvector if and only if the function $g_{y,\l} = \sum_{k=1}^{r-t+1} v_k f_{y,k+t-1} \in V_y$ is an eigenfunction of $A$ corresponding to the eigenvalue $\l$. The function $g_{y,\l}$ vanishes on $B(n,t-1) = \cup_{i=0}^{t-1} S(n,i)$. Furthermore, the restriction of $g_{y,\l}$ to each of the Hamming spheres $S(n,i)$, $t \le i \le r$ lies in the $t^{th}$ eigenspace $V_t$ of $S(n,i)$. In addition, it is easy to see that the tri-diagonal structure of $A_y$ implies that any non-zero eigenvector $v$ of $A_y$ has a non-zero first coordinate. Hence the restriction of the function $g_{y,\l}$ to $S(n,t)$ is a (non-zero) zonal spherical function around $y$ in the $t^{th}$ eigenspace $V_t$ of $S(n,t)$.

\noi We can now complete the proof of Theorem~\ref{thm:main}. Let $\l \in 2 \cdot {\cal R}(n-2t,r-t+1) - (n-2t)$. For each $y \in S(n,t)$, the space $V_y$ contains a $1$-dimensional eigenspace $U_y(\l)$ of $A$ corresponding to eigenvalue $\l$. Let $W_t(\l) = \<U_y(\l)\>_{|y| = t}$. Then $W_t(\l)$ is an eigenspace of $A$  corresponding to the eigenvalue $\l$. Note that if $f \in W_{t}(\l)$, then $f$ vanishes on $B(n,t-1) = \cup_{i=0}^{t-1} S(n,i)$. Furthermore, the restriction of $f$ to each of the Hamming spheres $S(n,i)$, $t \le i \le r$ lies in the $t^{th}$ eigenspace $V_t$ of $S(n,i)$.

\noi We claim that the of dimension of $W_t(\l)$ is at least ${n \choose t} - {n \choose {t-1}}$. To see that, note that its restriction to $S(n,t)$ is given by the span of zonal spherical functions in $V_t$ around all points $y \in S(n,t)$. By Lemma~\ref{lem:spherical zonal}, this span is the whole eigenspace $V_t$, whose dimension (see Section~\ref{subsubsec:doubly transitive}) is ${n \choose t} - {n \choose {t-1}}$.

\noi Next, let $t_1 < t_2$ with $\l \in 2 \cdot {\cal R}\(n-2t_1,r-t_1+1\) - \(n-2 t_1\)$ and $\l \in 2 \cdot {\cal R}\(n-2t_2,r-t_2+1\) - \(n-2 t_2\)$. We claim that the spaces $W_{t_1}(\l)$ and $W_{t_2}(\l)$ are orthogonal to each other. In fact, let $f \in W_{t_1}(\l)$ and $g \in W_{t_2}(\l)$. For $0 \le i \le r$, let $f_i, g_i$ be the restrictions of $f$ and $g$ to $S(n,i)$. By the discussion above, $f_i$ and $g_i$ belong to distinct eigenspaces $V_{t_1}$ and $V_{t_2}$ of $S(n,i)$ and hence are orthogonal to each other. It follows that $f \perp g$.

\noi Next, let $t_1 < t_2$ with $\l \in \Lambda_{t_1}$ and $\l \in \Lambda_{t_2}$. We claim that the spaces $W_{t_1}(\l)$ and $W_{t_2}(\l)$ are orthogonal to each other. In fact, let $f \in W_{t_1}(\l)$ and $g \in W_{t_2}(\l)$. For $0 \le i \le r_2$, let $f_i, g_i$ be the restrictions of $f$ and $g$ to $S(n,i)$. By the discussion above, $f_i$ and $g_i$ belong to distinct eigenspaces $V_{t_1}$ and $V_{t_2}$ of $S(n,i)$ and hence are orthogonal to each other. It follows that $f \perp g$.

\noi Hence for a real number $\l$ in $\cup_{t=0}^{r} \Lambda_t$, the multiplicity of $\l$ as an eigenvalue of $A$ is at least $\sum_{t: \l \in \Lambda_t} \({n \choose t} - {n \choose {t-1}}\)$.

\noi Since
\[
\sum_{\l \in \cup_{t=0}^{r} \Lambda_t} ~~\sum_{t: \l \in \Lambda_t} \({n \choose t} - {n \choose {t-1}}\) ~=~
\sum_{t=0}^r \({n \choose t} - {n \choose {t-1}}\) \cdot |\Lambda_t| ~=~
\]
\[
\sum_{t=0}^{r} \({n \choose t} - {n \choose {t-1}}\) \cdot (r-t+1) ~=~ \sum_{t=0}^r {n \choose t},
\]
and since the dimension of $A$ is $\sum_{t=0}^r {n \choose t}$, this shows that the multiplicity of each $\l$ in $\cup_{t=0}^r \Lambda_t$ is precisely $\sum_{t: \l \in \Lambda_t} \({n \choose t} - {n \choose {t-1}}\)$, and that $A$ has no additional eigenvalues.

\noi In particular, it follows that the dimension of each $W_t(\l)$ is ${n \choose t} - {n \choose {t-1}}$. Since, as discussed above, this is also the dimension of the restriction of this space to $S(n,t)$, any function $f \in W_{t}(\l)$ is determined by its restriction to $S(n,t)$.

\noi This completes the proof of the theorem.

\eprf

\subsection{Proof of Corollary~\ref{cor:max eigenvalue}}
We first observe that since $A$ is the adjacency matrix of a connected graph, by the Perron-Frobenius theorem, its maximal eigenvalue has multiplicity $1$, and the corresponding eigenfunction is strictly positive. 

\noi Our next observation is that for any $1 \le i \le n/2$, the zeroth eigenspace $V_0$ of $S(n,i)$ is the space of constant functions. Hence, for any $1 \le t \le i$, the functions contained is in the $t^{th}$ eigenspace $V_t$ of $S(n,i)$ sum to $0$, and hence cannot be strictly positive. Therefore, the second claim of Theorem~\ref{thm:main} implies that the space of maximal eigenfunctions of $A$ is  necessarily of the form $W_0(\l)$, that is, it is constructed by choosing the value of the parameter $t$ to be $0$ in the theorem. Let $\l$ be the maximal eigenvalue. Then the corresponding eigenspace $W_0(\l)$ consists, by the second claim of the theorem, of functions whose restrictions to the Hamming spheres $S(n,i)$ lie in the zeroth eigenspace of the corresponding sphere, and hence are constant on this sphere. This means that $W_0(\l)$ contains spherical functions on $B$, and hence the maximal eigenfunction is a positive spherical function. 

\noi Finally note that, by the first claim of Theorem~\ref{thm:main}, the eigenvalues of $A$ which are obtained by choosing $t=0$ belong to the set $\Lambda_0 = 2 \cdot {\cal R}(n,r+1) - n$. The maximal eigenvalue is the largest number in this set, and hence it is given by $\l = 2y - n$, where $y$ is the maximal root of $K^{(n)}_{r+1}$. Since the roots of $K^{(n)}_{r+1}$ are symmetric around $n/2$ (\cite{Szego}), it can be also written as $\l = n -2x$, where $y$ is the minimal root of $K^{(n)}_{r+1}$. This completes the proof of the Corollary.

\eprf

\subsection{Proof of Corollary~\ref{cor:max-ev}}
Let $1 \le s \le 2^{n-1}$ and let $r = n H^{-1}\(\frac{\log_2(s)}{n}\)$. Let $t = \lfloor r \rfloor$. We may assume that $r$ and hence $t$ are sufficiently large. Let $B = B(n,t)$ be the Hamming ball of radius $t$ around $0$. Then $|B| \le 2^{H\(\frac{t}{n}\) \cdot n} \le 2^{H\(\frac{r}{n}\) \cdot n} = s$. For the first inequality see e.g., Theorem~1.4.5. in \cite{van Lint}. 

\noi Let $\l(B)$ denote the maximal eigenvalue of the subgraph of the Hamming cube induced by the vertices in $B$. Since $\Lambda(s)$ is clearly non-decreasing (in fact, increasing) in $s$, Corollary~\ref{cor:max eigenvalue} implies that $\Lambda(s) \ge \l(B) = n - 2x_{t+1}$, where $x$ is the first root of the Krawtchouk polynomial $K^{(n)}_{t+1}$. 

\noi Combining bounds (5.36) and (5.44) in \cite{Lev-chapter}, we obtain the following bound on $x_{t+1}$, for $t \le \frac n2$:
\[
x_{t+1} ~\le~ \frac n2 - \sqrt{(t+1)(n-t+1)} + O_{t \rarrow \infty}\(t^{-1/6} \cdot \sqrt{n}\).
\]   
Hence
\[
\Lambda(s) ~\ge~ 2\sqrt{(t+1)(n-t+1)} - O_{t \rarrow \infty}\(t^{-1/6} \cdot \sqrt{n}\) ~\ge~ 2\sqrt{r \(n-r\)} - O_{r \rarrow \infty}\(r^{-1/6} \cdot \sqrt{n}\).
\]
This completes the proof of the first part of the corollary. 

\noi We pass to the second claim of the corollary. Let $s = 2^{n - o(n)}$. Observe that by the first claim of the corollary, $\Delta(s) \le \ln\(\frac{2^n}{s}\) \cdot \(1 + o_n(1)\) + O\(n^{1/3}\)$. This, however, fails to prove the second claim if the error term $O\(n^{1/3}\)$ is larger than the main term $\ln\(\frac{2^n}{s}\)$. 

\noi As above let $r = n H^{-1}\(\frac{\log_2(s)}{n}\)$, and let $t = \lfloor r \rfloor$. Then $\Lambda(s) \ge n - 2x_{t+1}$, and hence $\Delta(s) \le 2x_{t+1}$. We will use properties of the Krawtchouk polynomials to tighten the bound on $x_{t+1}$ and hence on $\Delta(s)$ in this case. 

\noi Recall that we denote by $K_k = K^{(n)}_k$ the $k^{th}$ Krawtchouk polynomial on $\{0,1\}^n$, as defined in Section~\ref{subsubsec:doubly transitive}. Krawtchouk polynomials $\{K_k\}_{k=0}^n$ form a family of orthogonal polynomials (with respect to the binomial measure on $[0,n]$) and hence all of their roots are real and distinct and lie in the interval $(0,n)$. Furthermore, the interval between any two consecutive roots of $K_k$ contains an integer point (\cite{Szego}). Let $x_k$ denote the first root of $K_k$. Then the sequence $\{x_k\}_{k=1}^n$ in $(0,n)$ is strictly decreasing (see e.g., \cite{Lev-chapter}). One additional property that we need is the {\it reciprocity} of Krawtchouk polynomials: For any $0 \le i, j \le n$ holds ${n \choose j} K_i(j) = {n \choose i} K_j(i)$ (see e.g., \cite{mrrw}). It follows that if $x_i \le t+1$, for some $0 \le i \le n$, then $x_{t+1} \le i$. In fact, note that $K_i$ is non-positive in the interval between $x_i$ and its second root. This interval contains an 
integer point $j \le t+1$. This implies that $K_i(j) \le 0$ and hence, by reciprocity $K_j(i) \le 0$. Hence $x_{t+1} \le x_j \le i$. Establishing a suitable upper bound on $i$ will lead to the required upper bound on $x_{t+1}$.    

\noi Let $\tau = 1 - \frac{\log_2(s)}{n} = \frac{\log_2\(2^n/s\)}{n}$. Then, by assumption, $\tau = o_n(1)$. We have $t+1 \ge r = \(\frac12 - \sqrt{\frac{\ln(2) \cdot \tau}{2}} + o_n(\tau)\) \cdot n$ where, as above, the second estimate is a simple consequence of the properties of the binary entropy function. Next, by the above bound from \cite{Lev-chapter}, we have $x_i \le \frac{n}{2} - \sqrt{i(n-i+1)} + O\(\sqrt{n}\)$. Combining these two estimates, it is easy to see that there exists a function $\e(s,n)$ which goes to $0$ if $n \rarrow \infty$ and $\frac{2^n}{s} \rarrow \infty$ so that if $i \ge \big(1+\e(s,n)\big) \cdot \frac{\ln\(\frac{2^n}{s}\)}{2}$ then $x_i \le t+1$. It follows that 
\[
\Delta(s) ~\le~ 2i ~\le~ \ln\(\frac{2^n}{s}\) \cdot \(1 + o_{s/2^n \rarrow 0}(1)\),
\]
proving the second claim of the corollary. 

\eprf

\end{document}